# A Multi-step Piecewise Linear Approximation Based Solution for Load Pick-up Problem in Electrical Distribution System

Jingyang Yun, Yun Zhou*, Weidong Hu, Peichao Zhang, Zheng Yan, Donghan Feng

School of Electronic Information and Electrical Engineering,
Shanghai Jiao Tong University, Shanghai 200240, China

*Abstract*—The load pick-up (LPP) problem searches the optimal configuration of the electrical distribution system (EDS), aiming to minimize the power loss or provide maximum power to the load ends. The piecewise linearization (PWL) approximation method can be used to tackle the nonlinearity and nonconvexity in network power flow (PF) constraints, and transform the LPP model into a mixed-integer linear programming model (LPP-MILP model). However, for the PWL approximation based PF constraints, big linear approximation errors will affect the accuracy and feasibility of the LPP-MILP model's solving results. And the long modeling and solving time of the direct solution procedure of the LPP-MILP model may affect the applicability of the LPP optimization scheme. This paper proposes a multi-step PWL approximation based solution for the LPP problem in the EDS. In the proposed multi-step solution procedure, the variable upper bounds in the PWL approximation functions are dynamically renewed to reduce the approximation errors effectively. And the multi-step solution procedure can significantly decrease the modeling and solving time of the LPP-MILP model, which ensure the applicability of the LPP optimization scheme. For the two main application schemes for the LPP problem (i.e. network optimization reconfiguration and service restoration), the effectiveness of the proposed method is demonstrated via case studies using a real 13-bus EDS and a real 1066-bus EDS.

## NOMENCLATURE

*Parameters*

| | |
|---|---|
| $\Lambda$ | Number of discretization in the PWL function |
| $g$ | Iteration number of the multi-step solution procedure of the LPP-MILP model |
| $\phi_{y,\lambda}$ | Slope parameter of the $\lambda^{\text{th}}$ segment in the PWL function of the square value of $y$ |
| $M$ | Sufficiently big constant number |
| $S_{\text{bus}}, S_{\text{feeder}}$ | Set of buses/feeders |
| $S_{\text{bus}}^{\text{DG}}$ | Set of buses with distributed generation |
| $r_{ij}, x_{ij}, z_{ij}$ | Resistance/reactance/impedance of feeder $ij$ |
| $V_{\text{bus}}^{\text{norm}}$ | Nominal bus voltage magnitude |
| $V_{\text{bus}}^{\text{min}}, V_{\text{bus}}^{\text{max}}$ | Lower/upper bound of bus voltage magnitude |

*Corresponding author.
E-mail addresses: yunjingyang94@sjtu.edu.cn (Jingyang Yun), yun.zhou@sjtu.edu.cn (Yun Zhou), freedomster@sjtu.edu.cn (Weidong Hu), pczhang@sjtu.edu.cn (Peichao Zhang), yanz@sjtu.edu.cn (Zheng Yan), seed@sjtu.edu.cn (Donghan Feng)

| | |
|---|---|
| $I_{ij}^{\max}$ | Upper bound of current magnitude of feeder $ij$ |
| $P_{\mathrm{G},i}^{\min}, Q_{\mathrm{G},i}^{\min}$ | Lower bound of active/reactive generation power of bus $i$ |
| $P_{\mathrm{G},i}^{\max}, Q_{\mathrm{G},i}^{\max}$ | Upper bound of active/reactive generation power of bus $i$ |
| $P_{ij}^{\max}, Q_{ij}^{\max}$ | Upper bound of active/reactive power of feeder $ij$ |
| $P_{\mathrm{L},i}^{\mathrm{par}}, Q_{\mathrm{L},i}^{\mathrm{par}}$ | Active/reactive load parameter of bus $i$ |

*Functions*

| | |
|---|---|
| $f(y, \bar{y}, \Lambda)$ | PWL approximation function of the square value of $y$ |

*Binary variables*

| | |
|---|---|
| $v_i, w_{ij}$ | State variable of bus $i$ /feeder $ij$ |

*Continuous variables*

| | |
|---|---|
| $\Delta_{y,\lambda}, y^+, y^-$ | Auxiliary variables used in the PWL function |
| $V_i^{\mathrm{sqr}}$ | Voltage magnitude square variable of bus $i$ |
| $I_{ij}^{\mathrm{sqr}}$ | Current magnitude square variable of feeder $ij$ |
| $P_{\mathrm{G},i}, Q_{\mathrm{G},i}$ | Active/reactive generation power of bus $i$ |
| $P_{\mathrm{L},i}, Q_{\mathrm{L},i}$ | Active/reactive restored load of bus $i$ |
| $P_{ij}, Q_{ij}$ | Active/reactive power of feeder $ij$ |

*Others*

| | |
|---|---|
| $S_{\mathrm{feeder}}^{\mathrm{use}}$ | Set of feeders in use state |
| $E_{\mathrm{p},ij}, E_{\mathrm{q},ij}$ | Error index in approximation of active/reactive power square variable of feeder $ij$ |
| $E_{\mathrm{p}}^{\mathrm{m}}, E_{\mathrm{q}}^{\mathrm{m}}$ | Mean error index in approximations of active/reactive power square variables of feeders |

## 1 INTRODUCTION

The load pick-up (LPP) problem in the electrical distribution system (EDS) searches the optimal configuration of the EDS, aiming to minimize the power loss or provide as much power as possible to load ends. The LPP problem is the core sub-problem of both network reconfiguration optimization and service restoration (reconfiguration after fault or blackout) optimization of the EDS. One fundamental model of the LPP problem for network reconfiguration is presented in [1]. And based on [1], the model of the LPP problem for service restoration optimization is proposed in [2].

The LPP problem belongs to the category of mixed-integer problem, since operational status of feeder or bus is defined strictly as 0 or 1. Along with other continuous variables, objective function and constraints, the LPP problem can be formulated as a mixed-integer nonlinear programming (MINLP) problem, which is a non-deterministic polynomial-time hard (NP-hard) and non-convex programming problem. Referring to the models in [3], network PF constraints are the essential part of the constraints in the LPP problem model (LPP model in short). And the network PF constraints make up the majority of nonlinear constraints in the LPP problem.

To tackle the nonlinearity and nonconvexity in network PF constraints, based on the DistFlow equations [4], three different formulations of PF constraints (i.e. DC PF constraints [5], second-order cone programming (SOCP) based PF constraints [6], and piecewise linearization (PWL) approximation based PF constraints [7]) are constructed and used in the LPP model. The DC PF constraints cannot account for the magnitudes of bus voltage and the distribution of reactive power in the EDS, which may generate an impractical solution. For the SOCP based PF constraints, the single quadratic equality constraint in the PF constraints is relaxed to second-order conic constraint, and the LPP model can be transformed into a SOCP model [6]. And for the SOCP based PF constraints, the sufficient condition of the SOCP relaxation should be considered [8]. The SOCP relaxation optimality preconditions may hinder the applicability of the SOCP based PF constraints for some EDS optimization schemes. Furthermore, the duality gap for the SOCP model may be non-zero, which will hinder the recovery of the alternative current (AC) feasibility from the SOCP solutions [9].

To further reduce the complexity in PF constraints, quadratic variables in the PF constraints are approximated by the PWL functions, and constitute the PWL based power flow (PF) constraints [7]. The PWL approximation based PF constraints have been widely used in different system optimization schemes [10]-[14]. With the PWL approximation based PF constraints, the LPP model can be transformed into a mixed-integer linear programming (MILP) model, and can be solved by commercial solvers (e.g. Gurobi and CPLEX) effectively. Moreover, due to the good numerical performance and convergence, the reasonable and acceptable modeling accuracy reduction of linear programming do not limit the industrial application of MILP in some cases [12].

For the PWL approximation based PF constraints, the approximation errors in the PWL approximation of the quadratic variables in the PF constraints should be researched specially. Previous literatures hardly concentrate on these topics. If the parameter of number of discretization $\Lambda$ in the PWL function is set to a small value, big approximation errors will affect the accuracy and feasibility of the MILP model's solving results. If the number of discretization $\Lambda$ is set to a big value, the corresponding MILP model's variables and constraints meet dramatically increase. And the long modeling and solving time of the direct solution procedure of the MILP model may affect the applicability of the optimization scheme.

In this paper, to improve the performance of the PWL approximation method and enhance the applicability of the LPP optimization scheme, a multi-step PWL approximation based solution for the LPP problem in the EDS is proposed. Apart from the network power flow constraints, other essential constraints of the LPP problem are presented, and make up the entire LPP model. The PWL approximation method is introduced to linearize the network power flow constraints, and transform the LPP model into a MILP model (LPP-MILP model). And a multi-step solution procedure of the LPP-MILP model is established. The LPP-MILP model can be solved taking advantage of warm starts from the solution of the previous iteration, where the variable upper bounds in the PWL approximation functions to segment the active/reactive power of feeders are renewed iteratively as the LPP-MILP model is solved. Error indices are constructed to evaluate the accuracy of the proposed multi-step PWL approximation based solution for the LPP problem in the EDS.

For the rest of the paper, the mathematical formulas of the LPP model are presented in Section 2. The detailed multi-step PWL approximation based solution for the LPP model are analyzed in Section 3. And the solvability, robustness and convergence of the multi-step solution method are analyzed and proved. The LPP-MILP model complexity assessment and error indices to evaluate the accuracy of the proposed method are presented in Section 4. Case studies to demonstrate the performance of the proposed method are included in Section 5. Eventually, Section 6 concludes the paper.

The main contributions of this research are as follows:

1) A general LPP optimization model for the two main optimization schemes of the EDS (i.e. network reconfiguration and service restoration) is established;

2) An effective multi-step PWL based solution for the LPP model is proposed with the advantages in solvability, robustness and convergence. And for large EDS, the applicability of the multi-step method can also be ensured.

## 2 MATHEMATICAL FORMULAS OF THE LPP MODEL

### 2.1 Objective Function

For the two main different application schemes for the LPP problem in the EDS (i.e. LPP problem for network reconfiguration and LPP problem for service restoration), different objective functions are as follows:

$$\min \sum_{ij \in S_{\text{feeder}}} r_{ij} I_{ij}^{\text{sqr}} \tag{1}$$

$$\max \sum_{i \in S_{\text{bus}}} P_{L,i} \tag{2}$$

For the LPP problem for network reconfiguration, minimizing the power loss of the EDS is considered as the optimization goal in the LPP model, and the corresponding objective function is shown in (1) [15]. For the LPP problem for service restoration, shown in (2), maximizing the restored load is set as the objective function in the LPP model [2].

### 2.2 Constraints

In this part, constraints considered in the LPP model are presented in detail, including the network PF constraints, the radial topology constraints, and other essential constraints.

*1) Network PF constraints*

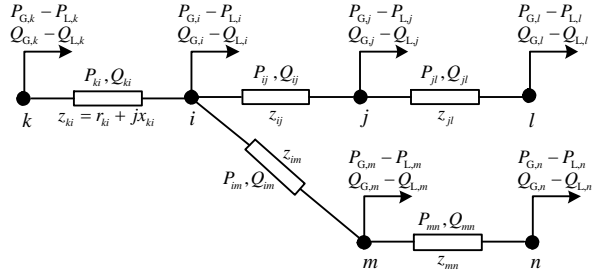

Fig. 1. Schematic diagram of a partial radial EDS network

Fig. 1 is the schematic diagram of a partial radial EDS network. For radial EDS network, the widely used DistFlow network PF constraints can be described as below [3]:

$$\sum_{ki \in S_{\text{feeder}}} P_{ki} - \sum_{ij \in S_{\text{feeder}}} (P_{ij} + r_{ij} I_{ij}^{\text{sqr}}) + P_{G,i} - P_{L,i} = 0 \quad i \in S_{\text{bus}} \tag{3}$$

$$\sum_{ki \in S_{\text{feeder}}} Q_{ki} - \sum_{ij \in S_{\text{feeder}}} (Q_{ij} + x_{ij} I_{ij}^{\text{sqr}}) + Q_{G,i} - Q_{L,i} = 0 \quad i \in S_{\text{bus}} \tag{4}$$

$$V_i^{\text{sqr}} - V_j^{\text{sqr}} = 2(P_{ij} r_{ij} + Q_{ij} x_{ij}) + (r_{ij}^2 + x_{ij}^2) I_{ij}^{\text{sqr}} \quad ij \in S_{\text{feeder}} \tag{5}$$

$$V_j^{\text{sqr}} I_{ij}^{\text{sqr}} = P_{ij}^2 + Q_{ij}^2 \quad ij \in S_{\text{feeder}} \tag{6}$$

$$(V_{\text{bus}}^{\min})^2 \le V_i^{\text{sqr}} \le (V_{\text{bus}}^{\max})^2 \quad i \in S_{\text{bus}} \tag{7}$$

$$I_{ij}^{\text{sqr}} \le (I_{ij}^{\max})^2 \quad ij \in S_{\text{feeder}} \tag{8}$$

The traditional equations of active and reactive power balance for each bus are represented as (3) and (4). The voltage difference across the feeder $ij$ can be obtained from (5). The relationship among the bus voltage magnitude, feeder current magnitude, and feeder active/reactive power is illustrated in (6), which is the only nonlinear constraint in the DistFlow network PF constraints. Limits for bus voltage magnitudes and feeder current magnitudes are shown in (7) and (8) respectively. Note that only quadratic terms of the magnitudes of bus voltage and feeder current (i.e. $V_i^{\text{sqr}}$ and $I_{ij}^{\text{sqr}}$) are used in (3)-(8), hence they shall not be treated as quadratic variables.

$$2(P_{ij} r_{ij} + Q_{ij} x_{ij}) - (r_{ij}^2 + x_{ij}^2) I_{ij}^{\text{sqr}} - (1 - w_{ij}) M \le V_i^{\text{sqr}} - V_j^{\text{sqr}} \le 2(P_{ij} r_{ij} + Q_{ij} x_{ij}) - (r_{ij}^2 + x_{ij}^2) I_{ij}^{\text{sqr}} + (1 - w_{ij}) M \quad ij \in S_{\text{feeder}} \tag{9}$$

$$P_{ij}^2 + Q_{ij}^2 - (1 - w_{ij}) M \le V_j^{\text{sqr}} I_{ij}^{\text{sqr}} \le P_{ij}^2 + Q_{ij}^2 + (1 - w_{ij}) M \quad ij \in S_{\text{feeder}} \tag{10}$$

$$(V_{\text{bus}}^{\min})^2 - (1 - v_i) M \le V_i^{\text{sqr}} \le (V_{\text{bus}}^{\max})^2 + (1 - v_i) M \quad i \in S_{\text{bus}} \tag{11}$$

$$I_{ij}^{\text{sqr}} \le (I_{ij}^{\max})^2 + (1 - w_{ij}) M \quad ij \in S_{\text{feeder}} \tag{12}$$

By introducing the state variables of buses and feeders, constraints (5)-(8) can be reformulated in (9)-(12). The state variables $v_i$ and $w_{ij}$ assumes a value of 1 if the bus or feeder is in the use state, and 0 if the bus or feeder is in the non-use state. If the use state of any bus or feeder is non-optional, the corresponding state variable $v_i$ or $w_{ij}$ should be set to a fixed value. Referring to [4], by introducing a sufficiently big constant number $M$, constraints in (9)-(12) can be relaxed if the corresponding bus or feeder is in the non-use state. In the LPP model, (3)-(4) and (9)-(12) make up the integrated network PF constraints.

*2) Radial topology constraints and other constraints*

Considering $N_s$ islands in the network, radial topology constraint for the EDS is shown in (13).

$$\sum_{ij \in S_{\text{feeder}}} w_{ij} = \sum_{i \in S_{\text{bus}}} v_i - N_s \tag{13}$$

$$v_i + v_j \geq 2w_{ij} \quad ij \in S_{\text{feeder}} \tag{14}$$

$$v_i P_{G,i}^{\min} \leq P_{G,i} \leq v_i P_{G,i}^{\max} \quad i \in S_{\text{bus}} \tag{15}$$

$$v_i Q_{G,i}^{\min} \leq Q_{G,i} \leq v_i Q_{G,i}^{\max} \quad i \in S_{\text{bus}} \tag{16}$$

$$v_i P_{L,i}^{\text{par}} \leq P_{L,i} \leq v_i P_{L,i}^{\text{par}} \quad i \in S_{\text{bus}} \tag{17}$$

$$v_i Q_{L,i}^{\text{par}} \leq Q_{L,i} \leq v_i Q_{L,i}^{\text{par}} \quad i \in S_{\text{bus}} \tag{18}$$

$$-w_{ij} M \leq P_{ij} \leq w_{ij} M \quad ij \in S_{\text{feeder}} \tag{19}$$

$$-w_{ij} M \leq Q_{ij} \leq w_{ij} M \quad ij \in S_{\text{feeder}} \tag{20}$$

Other essential constraints in the LPP model are listed in (14)-(20). The relationship between bus and feeder state variables are indicated in (14). When the feeder $ij$ is in use state, the corresponding endpoints of this feeder are accordingly in use state (i.e. $v_i$=1 and $v_j$=1 ). Restricted by the state variables, the limits of $P_{G,i}$, $Q_{G,i}$, $P_{L,i}$, $Q_{L,i}$, $P_{ij}$, and $Q_{ij}$ are shown in (15)-(20) respectively.

## 3 MULTI-STEP PWL APPROXIMATION BASED SOLUTION FOR THE LPP MODEL

### 3.1 PWL approximation of network PF constraints

In the LPP model presented in Section 2, only constraint (6) or (10) (considering the state variables) contains quadratic variables and thus is nonlinear in the network PF constraints. As a standard procedure in electrical system optimization problems, the PWL approximation method is used to process the quadratic variables in (6) or (10), and to linearize the nonlinear network power flow constraint [7]. To make the presentation clearer, constraint (6) is used to illustrate the PWL approximation procedure:

$$V_j^{\text{sqr}} I_{ij}^{\text{sqr}} \approx (V_{\text{bus}}^{\text{norm}})^2 I_{ij}^{\text{sqr}} \quad ij \in S_{\text{feeder}} \tag{21}$$

$$(V_{\text{bus}}^{\text{norm}})^2 I_{ij}^{\text{sqr}} = P_{ij}^2 + Q_{ij}^2 \quad ij \in S_{\text{feeder}} \tag{22}$$

$$P_{ij}^2 \approx f(P_{ij}, P_{ij}^{\max}, \Lambda) \quad ij \in S_{\text{feeder}} \tag{23}$$

$$Q_{ij}^2 \approx f(Q_{ij}, Q_{ij}^{\max}, \Lambda) \quad ij \in S_{\text{feeder}} \tag{24}$$

$$(V_{\text{bus}}^{\text{norm}})^2 I_{ij}^{\text{sqr}} = f(P_{ij}, P_{ij}^{\max}, \Lambda) + f(Q_{ij}, Q_{ij}^{\max}, \Lambda) \quad ij \in S_{\text{feeder}} \tag{25}$$

Firstly, the $V_j^{\text{sqr}}$ in the left side of (6) is approximated by $(V_{\text{bus}}^{\text{norm}})^2$ (shown in (21)), and constraint (6) is reformulated as (22). Due to the narrow voltage magnitude interval [$V_{\text{bus}}^{\min}$, $V_{\text{bus}}^{\max}$] in the EDS, the simplification of (21) has a relatively low approximation error [7]. Secondly, the quadratic variables $P_{ij}^2$ and $Q_{ij}^2$ in (21) are approximated by PWL functions $f(y, \bar{y}, \Lambda)$ (shown in (23) and (24)), and then constraint (6) can be further expressed as (25). Similarly, considering the state variables of the EDS, constraint (10) can be formulated as (26).

$$f(P_{ij}, P_{ij}^{\max}, \Lambda) + f(Q_{ij}, Q_{ij}^{\max}, \Lambda) - (1 - w_{ij}) M \leq (V_j^{\text{norm}})^2 I_{ij}^{\text{sqr}} \leq f(P_{ij}, P_{ij}^{\max}, \Lambda) + f(Q_{ij}, Q_{ij}^{\max}, \Lambda) + (1 - w_{ij}) M \quad ij \in S_{\text{feeder}} \tag{26}$$

The PWL function $f(y,\bar{y},\Lambda)$ is defined in (27)-(32), and Fig. 2 shows how the PWL function approximates the quadratic curve. Considering the paper length, the detailed presentation of the PWL function can be referred to [7].

$$f(y,\bar{y},\Lambda) = \sum_{\lambda=1}^{\Lambda} \phi_{y,\lambda} \Delta_{y,\lambda} \tag{27}$$

$$y = y^+ - y^- \tag{28}$$

$$y^+ + y^- = \sum_{\lambda=1}^{\Lambda} \Delta_{y,\lambda} \tag{29}$$

$$0 \leq \Delta_{y,\lambda} \leq \bar{y}/\Lambda \quad \forall \lambda=1,2,...,\Lambda \tag{30}$$

$$\phi_{y,\lambda} = (2\lambda-1)\bar{y}/\Lambda \quad \forall \lambda=1,2,...,\Lambda \tag{31}$$

$$y^+, y^- \geq 0 \tag{32}$$

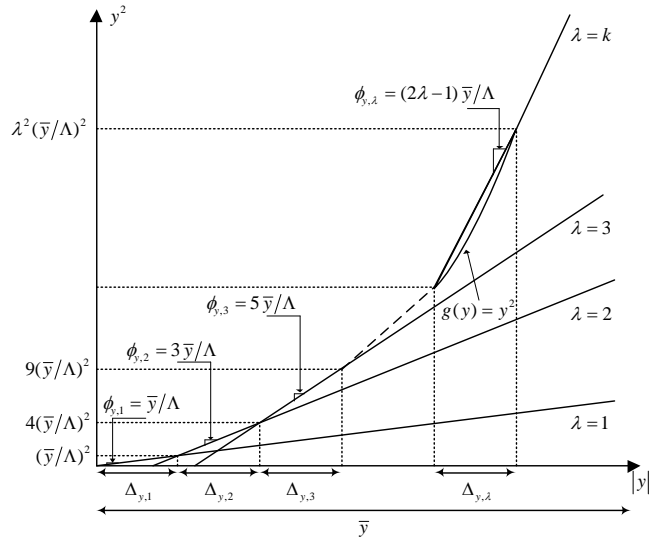

Fig. 2. Piecewise linear approximation function

By linearization of the network PF constraints, the LPP model is transformed into a MILP model with objective function in (1) or (2) and constraints in (3)-(4), (9), (11)-(20), and (26).

*3.2 Multi-step solution procedure of the LPP-MILP model*

The errors in the PWL approximation of the network PF constraints may affect the feasibility of the solving results of the LPP-MILP model. In general, when implementing the PWL approximation method, bigger value of $\Lambda$ contributes to better fitting of the quadratic curves of $P_{ij}^2$ and $Q_{ij}^2$, and therefore the approximation and the solution will be more accurate. However, the solution efficiency may be meanwhile sacrificed severely due to large scale variables generated by the PWL functions, which sharply increases the dimensions of the corresponding optimization problem.

$$P_{ij}^{\max} = V_{\text{bus}}^{\max} I_{ij}^{\max} \tag{33}$$

$$Q_{ij}^{\max} = V_{\text{bus}}^{\max} I_{ij}^{\max} \tag{34}$$

Except for the number of discretization $\Lambda$, the set value of upper limit (e.g. $P_{ij}^{\max}$ or $Q_{ij}^{\max}$ in (26)) in the PWL function is another important factor influencing the approximation accuracy of the quadratic variable. For given $V_{\text{bus}}^{\max}$ and $I_{ij}^{\max}$, $P_{ij}^{\max}$ and $Q_{ij}^{\max}$ can be calculated by (33) and (34) respectively. For a certain $I_{ij}^{\max}$ for a EDS (e.g. $I_{ij}^{\max}$=200 A for a 44-bus test system in [7]), the actual active/reactive power through feeder $P_{ij}$ and $Q_{ij}$ may be far less than $P_{ij}^{max}$ and $Q_{ij}^{max}$ in some cases, especially in the LPP problem due to the limited generating capacity in the EDS. The big gaps between the actual values and the upper limit set values will affect the accuracy of the PWL approximation method heavily.

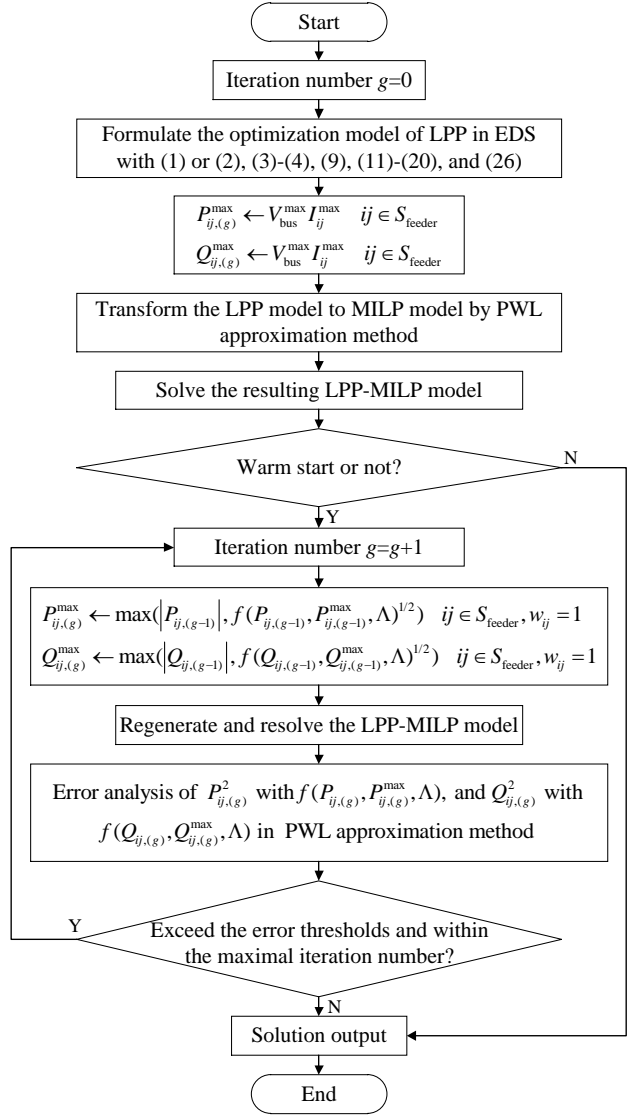

Fig. 3. Flow chart of the multi-step solution procedure of the LPP-MILP model

In the solution of the LPP-MILP model, for a small set value of $\Lambda$, if upper limit set values can be rated closer to the actual values, the accuracy of the PWL approximation method can be improved, meanwhile the computational complexity is not increased. Taking advantages of warm starts from the solution of the previous iteration, a multi-step solution procedure of the LPP-MILP model is proposed. The flow chart of the multi-step solution procedure of the LPP-MILP model is shown in Fig. 3. In this flowchart, the symbols with subscript ($g$) are belong to the parameters and solving results of the LPP-MILP model under the $g$-th iteration solution step.

In the multi-step solution flow chart, the initial generated LPP-MILP model is solved firstly. And then the active/reactive power upper limit set values of feeders in use state ($w_{ij}=1$) are renewed by the results solved in the first step. The renewal principles of $P_{ij,(g)}^{max}$ and $Q_{ij,(g)}^{max}$ designed in Fig. 3 can ensure the solvability, and robustness, and convergence of the regenerated LPP-MILP model. After error analysis, if errors between $P^2_{ij,(g)}$ and $f(P_{ij,(g)}, P^{\max}_{ij,(g)}, \Lambda)$ or $Q^2_{ij,(g)}$ and $f(Q_{ij,(g)}, Q^{\max}_{ij,(g)}, \Lambda)$ exceed the error threshold and within the maximal iteration number, the iteration solution step can be reduplicated.

*3.3 Solvability, robustness, and convergence analysis*

The solvability, robustness, and convergence of the proposed multi-step solution procedure of the LPP-MILP model are analyzed in this part. For narrative convenience, subsequent analysis take the active power part in Fig. 3 for example. The reactive power part can be analyzed in a similar way.

*1) Solvability*

$$P^{\max}_{ij,(0)} = V^{\max}_{\text{bus}} I^{\max}_{ij} \tag{35}$$

In Fig. 3, in the initial iteration solution step ($g = 0$), for feeder $ij$, the feeder active power upper bound parameter $P^{\max}_{ij(0)}$ is set to $V^{\max}_{\text{bus}} I^{\max}_{ij}$ in (35), where $I^{\max}_{ij}$ is the upper bound of feeder current magnitude based on the maximal acceptable temperature of the conductors. And with other parameters, the solvability of the LPP-MILP model in the initial iteration solution step can be ensured.

*2) Robustness*

In the proposed multi-step solution procedure for the LPP-MILP model, the variable upper bounds of feeder power flows ($P^{\max}_{ij,(g)}$) in the PWL approximation functions are dynamically renewed. For the LPP-MILP model and solving results in the $g$-th step, only the feeder power flow $P_{ij,(g)}$ and its quadratic approximation term $f(P_{ij,(g)}, P^{\max}_{ij,(g)}, \Lambda)$ are changed.

$$|P_{ij,(g-1)}| \leq P^{\max}_{ij,(g)} \tag{36}$$

$$f(P_{ij,(g-1)}, P^{\max}_{ij,(g-1)}, \Lambda) \leq f(P^{\max}_{ij,(g)}, P^{\max}_{ij,(g)}, \Lambda) \tag{37}$$

To guarantee the robustness of the iteration solution step in the multi-step solution procedure, the solving results of the last iteration solution step should be within the solution space of the LPP-MILP model in the current iteration solution step. Referring to the constraints involved in the LPP-MILP model (i.e. (3)-(4), (9), (11)-(20), and (26)), conditions in (36) and (37) should be satisfied to guarantee the robustness.

$$f(P^{\max}_{ij,(g)}, P^{\max}_{ij,(g)}, \Lambda)^{1/2} = P^{\max}_{ij,(g)} \tag{38}$$

From the equations of the PWL function in (27)-(32), the equality relationship in (38) is correct. And with the renewal strategy of $P^{\max}_{ij,(g)}$ designed in Fig. 3, conditions in (36) and (37) can be demonstrated. Accordingly, the robustness of the iteration solution step in the multi-step solution procedure can be guaranteed.

*3) Convergence*

To analysis the convergence of the proposed method, two lemmas are introduced and proved as following:

*Lemma 1*:

$$|P_{ij,(g)}| \leq f(P_{ij,(g)}, P^{\max}_{ij,(g)}, \Lambda)^{1/2} \tag{39}$$

*Proof 1*: Fig. 2 shows how the PWL function approximates the quadratic curve $h(y) = y^2$. Referring to (27)-(32), due to the convexity of $h(y)$, all segments employed in the PWL function are above on $h(y)$ (shown in Fig. 2), and the Lemma 1 is correct obviously.

*Lemma 2*:

$$P^{\max}_{ij,(g)} \leq P^{\max}_{ij,(g-1)} \tag{40}$$

*Proof 2*: Combined with *Lemma 1* and the renewal strategy of $P^{\max}_{ij,(g)}$ in Fig. 3, (40) is equivalent to:

$$f(P_{ij,(g-1)}, P^{\max}_{ij,(g-1)}, \Lambda)^{1/2} \leq P^{\max}_{ij,(g-1)} \tag{41}$$

From the equations of the PWL function in (27)-(32), (41) is correct, and accordingly *Lemma 2* can be demonstrated.

$$\sup\{|f(P_{ij,(g)}, P_{ij,(g)}^{\max}, \Lambda) - P_{ij,(g)}^2|\} \le \sup\{|f(P_{ij,(g)}, P_{ij,(g-1)}^{\max}, \Lambda) - P_{ij,(g)}^2|\} \tag{42}$$

In (42), $P_{ij,(g)}$ belongs to the solving results of the feeder active power in the $g$-th iteration solution step. And referring to (23), $|f(P_{ij,(g)}, P_{ij,(g)}^{\max}, \Lambda) - P_{ij,(g)}^2|$ in (42) can be defined as the PWL approximation error item of $P_{ij,(g)}^2$. Considering all possible segmented state of $P_{ij,(g)}$ from (28)-(32), with Lemma 1 and 2, (42) can be demonstrated.

As the supremum of the approximation error item of $P_{ij,(g)}^2$ declines with the increasing of the iteration solution steps, combined with the previous robust analysis, the convergence of the proposed multi-step solution procedure can be guaranteed, while detailed convergence property will be further stressed in the case studies.

## 4 LPP-MILP MODEL COMPLEXITY ASSESSMENT AND ERROR INDICES CONSTRUCTION

To assess the complexity of the LPP-MILP model, the counts of total variables and constraints in the LPP-MILP model are analyzed and listed in Table I. As shown in Table I, the total variables and constraints of the LPP-MILP model are positively related to the number of discretization $\Lambda$ in the PWL function. For a large EDS with big $|S_{\text{feeder}}|$ and $|S_{\text{bus}}|$, if a relatively big value of $\Lambda$ is taken into use, the curse of dimensionality of the LPP-MILP model may lead to an impractical long modeling and solving time of the LPP-MILP model. But for the multi-step solution procedure in Fig. 3, the number of variables and constraints in the LPP-MILP model in each iterative solution step is invariable. As a small set value of $\Lambda$ can be used in the multi-step solution procedure of the LPP-MILP model, the curse of dimensionality problem can be avoided. And for a large EDS, the computing complexity of the LPP-MILP model can be reduced.

TABLE I
COUNTS OF VARIABLES AND CONSTRAINTS OF THE LPP-MILP MODEL

| Variables | | Constraints | |
|---|---|---|---|
| Binary variables | Continuous variables | Network PF constraints | Radial topology constraints and other constraints |
| $v_i, w_{ij}$ | $P_{G,i}, Q_{G,i}, P_{L,i}, Q_{L,i}, P_{ij}, Q_{ij},$ $P_{ij}^+, P_{ij}^-, Q_{ij}^+, Q_{ij}^-, \Delta_{P_{ij},\lambda}, \Delta_{Q_{ij},\lambda}$ $V_i^{\text{sqr}}, I_{ij}^{\text{sqr}}$ | (3)-(4), (9), (11)-(12), (26) | (13)-(20) |
| $|S_{\text{feeder}}| + |S_{\text{bus}}|$ | $(7+2\Lambda)|S_{\text{feeder}}| + 5|S_{\text{bus}}|$ | $(13+2\Lambda)|S_{\text{feeder}}| + 3|S_{\text{bus}}|$ | $3|S_{\text{feeder}}| + 4|S_{\text{bus}}| + 1$ |

The PWL approximation method of network PF constraints has been presented in Section 3.1. The approximations in (21), (23) and (24) make up the core steps of the PWL approximation method. As the simplification of (21) has a relatively low approximation error [7], the errors caused by approximations in (23) and (24) mainly determine the accuracy of the PWL approximation method in the LPP-MILP model.

$$S_{\text{feeder}}^{\text{use}} = \{ij | ij \in S_{\text{feeder}}^{\text{f}} \text{ and } w_{ij}=1\} \tag{43}$$

$$E_{p,ij} = \frac{|f(P_{ij}, P_{ij}^{\max}, \Lambda) - P_{ij}^2|}{P_{ij}^2} \times 100\% \quad ij \in S_{\text{feeder}}^{\text{use}} \tag{44}$$

$$E_{q,ij} = \frac{|f(Q_{ij}, Q_{ij}^{\max}, \Lambda) - Q_{ij}^2|}{Q_{ij}^2} \times 100\% \quad ij \in S_{\text{feeder}}^{\text{use}} \tag{45}$$

$$E_p^m = \frac{1}{|S_{\text{feeder}}^{\text{use}}|} \sum_{ij \in S_{\text{feeder}}^{\text{use}}} E_{p,ij} \tag{46}$$

$$E_q^m = \frac{1}{|S_{\text{feeder}}^{\text{use}}|} \sum_{ij \in S_{\text{feeder}}^{\text{use}}} E_{q,ij} \tag{47}$$

For the solving results of the LPP-MILP model, the set of feeders in use state ($w_{ij} = 1$) is defined firstly in (43). For the feeder $ij$, error indices of approximations in (23) and (24) ($E_{p,ij}$ and $E_{q,ij}$) are defined in (44) and (45) respectively. To estimate the global accuracy of the PWL approximation method in the LPP-MILP model, the mean error indices ($E_p^m$ and $E_q^m$) are defined in (46) and (47) respectively.

## 5 CASE STUDIES

In this section, the proposed multi-step solution procedure for the LPP problem is tested using the following two real electrical distribution systems.

a) A real 13-bus distribution system with 2 DGs in the University of Manchester, UK [16].

b) A real 1066-bus distribution system with 44 DGs in Shandong province of China [17].

For each EDS test case, the two main application schemes for the LPP problem (i.e. network reconfiguration and service restoration) are tested respectively. The case studies are performed on a 3.2-GHz dual-core PC with 16GB RAM. The LPP-MILP model and the multi-step solution procedure are implemented in MATLAB with YALMIP toolbox, and Gurobi 8.0.1 is the MILP solver.

### 5.1 Case I: the 13-bus EDS with 2 DGs

The University of Manchester owns a 6.6 kV distribution network fed from a 33/6.6 kV transformer connected to bus 13. Fig. 4 shows the topology of the 13-bus EDS with 2 DG in the University of Manchester, while the actual area in the university served by the electrical network can be referred to [16].

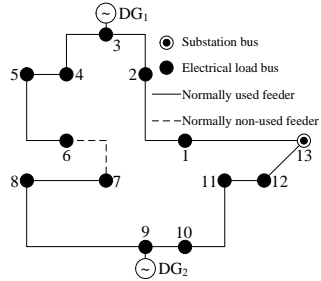

Fig. 4. The topology of the 13-bus EDS with 2 DGs in the University of Manchester, UK

TABLE II
DG PARAMETERS OF THE 13-BUS DISTRIBUTION SYSTEM

| No. | Name | $P_{G,i}^{min}$ (MW) | $P_{G,i}^{max}$ (MW) | $Q_{G,i}^{min}$ (Mvar) | $Q_{G,i}^{max}$ (Mvar) |
|---|---|---|---|---|---|
| 1 | DG$_1$ | 0 | 0.7497 | -0.2721 | 0.2721 |
| 2 | DG$_2$ | 0 | 0.8302 | -0.3014 | 0.3014 |

TABLE III
ELECTRICAL BUS LOAD PARAMETERS OF THE 13-BUS DISTRIBUTION SYSTEM

| Bus | $P_{L,i}^{par}$ (MW) | $Q_{L,i}^{par}$ (MW) | Bus | $P_{L,i}^{par}$ (MW) | $Q_{L,i}^{par}$ (MW) |
|---|---|---|---|---|---|
| 1 | 0.2420 | 0.0878 | 7 | 0.1614 | 0.0586 |
| 2 | 0.4775 | 0.1733 | 8 | 0.5180 | 0.1880 |
| 3 | 0.1068 | 0.0388 | 9 | 0.0684 | 0.0248 |
| 4 | 0.1772 | 0.0643 | 10 | 0.0953 | 0.0346 |
| 5 | 0.1242 | 0.0451 | 11 | 0.2137 | 0.0776 |
| 6 | 0.1920 | 0.0697 | 12 | 0.0588 | 0.0213 |

The system base voltage is 6.6 kV and the base capacity is 1 MVA. The DG and bus load parameters are derived from [16] and listed in Table II and III, and the detailed feeder parameters of the test system can be referred to [16]. The bus voltage and feeder current limits are set as follows: $V_{\text{bus}}^{\min}=0.95\ V_{\text{bus}}^{\text{norm}}$, $V_{\text{bus}}^{\max}=1.05\ V_{\text{bus}}^{\text{norm}}$, and $I_{ij}^{\max}=250$ A $ij \in S_{\text{feeder}}$. In Case I, the optimality gap for MILP solver is set as 0.01%.

*1) Network reconfiguration scheme*

The LPP-MILP model for the network reconfiguration scheme includes objective function in (1) and constraints in (3)-(4), (9), (11)-(20), and (26). Besides, all bus loads should be served in the network reconfiguration scheme.

TABLE IV
SOLVING RESULTS OF THE LPP-MILP MODEL BY THE DIRECT SOLUTION PROCEDURE UNDER DIFFERENT SET VALUES OF $\Lambda$ FOR THE NETWORK RECONFIGURATION SCHEME (CASE I)

| No. | Set values of $\Lambda$ | Modeling and solving time (s) | Objective function value (kW) | $E_{\text{p}}^{\text{m}}$ (%) | $E_{\text{q}}^{\text{m}}$ (%) |
|---|---|---|---|---|---|
| 1  | $\Lambda=10$  | 2.1221  | 1.0359 | 107.021399 | 369.205235 |
| 2  | $\Lambda=20$  | 3.6104  | 0.9946 | 46.276745  | 156.355773 |
| 3  | $\Lambda=30$  | 5.3291  | 0.9875 | 27.569278  | 93.206063  |
| 4  | $\Lambda=40$  | 7.3771  | 0.9853 | 18.347546  | 65.530572  |
| 5  | $\Lambda=50$  | 9.2774  | 0.9844 | 12.775294  | 50.042951  |
| 6  | $\Lambda=60$  | 11.5417 | 0.9835 | 9.276318   | 40.084584  |
| 7  | $\Lambda=70$  | 14.1038 | 0.9831 | 6.632350   | 33.051539  |
| 8  | $\Lambda=80$  | 16.9370 | 0.9829 | 4.817009   | 27.874825  |
| 9  | $\Lambda=90$  | 19.9414 | 0.9827 | 3.307448   | 23.213685  |
| 10 | $\Lambda=100$ | 23.1952 | 0.9825 | 2.121463   | 19.915755  |

For the network reconfiguration scheme of the 13-bus EDS, Table IV shows the solving results of the LPP-MILP model with the direct solution procedure (i.e. multi-step iteration procedure not selected in Fig. 3) under different set values of $\Lambda$. Referring to Table IV, bigger set value of $\Lambda$ can lead to lower error index values, and promote better accuracy of the PWL approximation method basically. The complexity assessment of the LPP-MILP model has been analyzed and listed in Table I, with the growth of the set value of $\Lambda$, the counts of variables and constraints of the LPP-MILP model meet dramatically increase. Accordingly, demonstrated in Table IV, the modeling and solving time of the LPP-MILP model increases significantly with the growth of the set value of $\Lambda$. Referring to [7] and select a more stringent set value of $\Lambda$, relaxing the error threshold and the maximum iteration number limit in Fig. 3, Table V shows the solving results of the LPP-MILP model by the multi-step solution procedure under different iteration number with fixed set value of $\Lambda$ ($\Lambda=10$). In Table V, the No. 1 scenario (iteration number=0) is just the solving results without the multi-step solution procedure in Table IV under $\Lambda=10$.

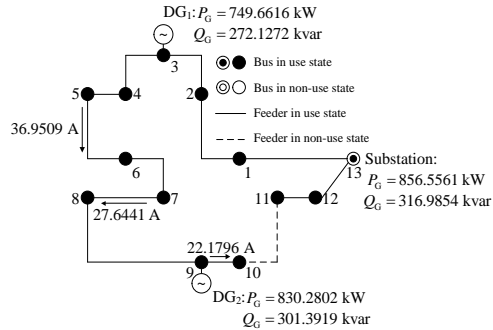

Fig. 5. Network reconfiguration results of the 13-bus EDS without the multi-step solution procedure ($\Lambda=10$, iteration number=0)

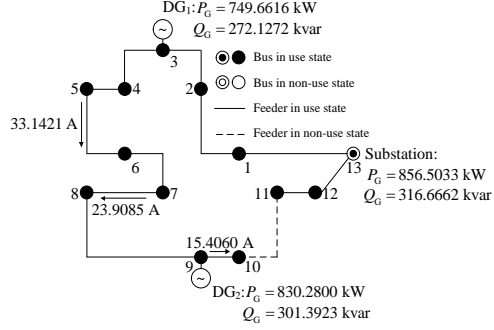

Fig. 6. Network reconfiguration results of the 13-bus EDS by the multi-step solution procedure ($\Lambda=10$, iteration number=1)

TABLE V
SOLVING RESULTS OF THE LPP-MILP MODEL BY THE MULTI-STEP SOLUTION PROCEDURE UNDER DIFFERENT ITERATION NUMBER ($\Lambda=10$) FOR THE NETWORK OPTIMIZATION RECONFIGURATION SCHEME (CASE I)

| No. | Iteration number | Accumulated modeling and solving time (s) | Objective function value (kW) | $E_p^m$ (%) | $E_q^m$ (%) |
|---|---|---|---|---|---|
| 1 | 0 | 2.1221 | 1.0359 | 107.021399 | 369.205235 |
| 2 | 1 | 2.3935 | 0.9830 | 0.294074 | 0.964232 |
| 3 | 2 | 2.6243 | 0.9821 | 0.014013 | 0.055010 |
| 4 | 3 | 2.8479 | 0.9820 | 0.000701 | 0.002719 |
| 5 | 4 | 3.1002 | 0.9820 | 0.000072 | 0.000096 |
| 6 | 5 | 3.3528 | 0.9820 | 0.000001 | 0.000057 |

The detailed network reconfiguration results of the 13-bus test system without the multi-step solution procedure (No. 1 scenario in Table V, $\Lambda=10$, iteration number=0) and by the multi-step solution procedure (No. 2 scenario in Table V, $\Lambda=10$, iteration number=1) are shown in Fig. 5 and 6 respectively. After one iteration process, the error index values have a significant decline, and the accuracy performance is better than the No.10 scenario ($\Lambda=100$) in Table IV with a much shorter modeling and solving time. And comparing the feeder current results of some representative feeders in Fig. 5 and 6, without the multi-step solution procedure (No. 1 scenario in Table V), except the deviation in objective function values, the heavily overestimation of the feeder currents in Fig. 5 due to big linear approximation errors will also affect the feasibility of the reconfiguration scheme solving results.

After two iteration processes (No. 3 scenario in Table V), the error indices are effectively decreased to relative small values ($E_p^m$, $E_q^m$ <0.1%). Solving results in Table V show that with the increasing of iteration numbers, the error index values reveal a trend of gradual decrease, which can demonstrate the effectiveness, and robustness of the proposed multi-step solution procedure of the LPP-MILP model. The good convergence of the multi-step solution procedure are also shown and demonstrated in Table V. Comparing the modeling and solving time in Table IV and V, the advantage in solution efficiency of the proposed methodology can be demonstrated.

*2) Service restoration optimization scheme*

The LPP-MILP model for the service restoration optimization scheme is formulated with objective function in (2) and constraints in (3)-(4), (9), (11)-(20), and (26). For the service restoration optimization scheme in Case I, the substation bus (bus 13) in Fig. 4 is set as disconnected with the upper transmission network.

For the solution procedure in Fig. 3, set the error thresholds as $E_p^m \leq 0.1\%$ and $E_q^m \leq 0.1\%$, and set the maximal iteration number as 5. Solving processes and results of the LPP-MILP model by the multi-step solution procedure are listed in Table VI ($\Lambda=10$). After 2 iteration processes, the error indices of the solving results are reduced within the error thresholds, and the multi-step solution procedure is completed. The corresponding service restoration optimization results are shown in Fig. 7.

TABLE VI
SOLVING RESULTS OF THE LPP-MILP MODEL BY THE MULTI-STEP SOLUTION PROCEDURE UNDER DIFFERENT ITERATION NUMBER ($\Lambda$ =10) IN THE SERVICE RESTORATION OPTIMIZATION SCHEME (CASE I)

| No. | Iteration number | Accumulated modeling and solving time (s) | Objective function value (kW) | $E_p^m$ (%) | $E_q^m$ (%) |
|---|---|---|---|---|---|
| 1 | 0 | 2.2668 | 1.4434 | 100.365949 | 361.380124 |
| 2 | 1 | 2.5224 | 1.4434 | 0.297782 | 0.297782 |
| 3 | 2 | 2.7798 | 1.4434 | 0.021539 | 0.028453 |

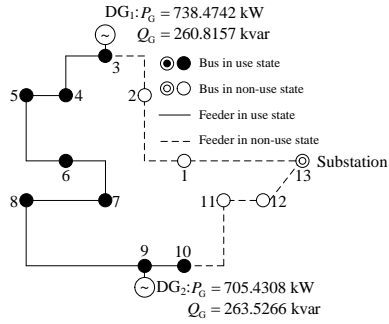

Fig. 7. Service restoration optimization results of the 13-bus EDS by the multi-step solution procedure ($\Lambda$ =10, iteration number=2)

TABLE VII
SOLVING RESULTS OF THE LPP-MILP MODEL BY THE DIRECT SOLUTION PROCEDURE UNDER DIFFERENT SET VALUES OF $\Lambda$ FOR THE SERVICE RESTORATION OPTIMIZATION SCHEME (CASE I)

| No. | Set values of $\Lambda$ | Modeling and solving time (s) | Objective function value (kW) | $E_p^m$ (%) | $E_q^m$ (%) |
|---|---|---|---|---|---|
| 1 | $\Lambda$ =100 | 23.2458 | 1.4434 | 5.071336 | 38.161921 |
| 2 | $\Lambda$ =500 | 491.7814 | 1.4434 | 0.099766 | 0.498889 |
| 3 | $\Lambda$ =900 | 1977.3171 | 1.4434 | 0.057309 | 0.090220 |

Referring to Table VII, for the direct solution procedure of the LPP-MILP model, to satisfy the same error precision (i.e. $E_p^m \leq 0.1\%$ and $E_q^m \leq 0.1\%$), the set value of $\Lambda$ need to increase to 900. Comparing the modeling and solving time in Table VI and VII, the advantage of the proposed multi-step solution procedure in solution efficiency can be demonstrated.

*5.2 Case II: the 1066-bus EDS with 44 DGs*

To stress the issue of application in relative large EDS, the proposed multi-step solution procedure for the LPP problem is tested on a real distribution system in Shandong province of China in Case II. The test system consists of 1066 buses and 44 DGs, and the detailed system data can be referred to [18]. As bus equivalent shunt admittance parameters included in the data of 1066-bus EDS, the active and reactive power balance constraints (i.e. (3) and (4) in Section 2.2) in the LPP-MILP model need to be improved. The supplemented power balance constraints can be referred to [10].

For the 1066-bus EDS, the bus voltage, feeder current limits, and number of discretization in the PWL function are set as follows: $V_{bus}^{min}$ =0.95 $V_{bus}^{norm}$, $V_{bus}^{max}$ =1.05 $V_{bus}^{norm}$, $I_{ij}^{max}$ =500 A $ij \in S_{feeder}$, and $\Lambda$ =10. For the solution procedure in Fig. 3, set the error thresholds as $E_p^m \leq 0.1\%$ and $E_q^m \leq 0.1\%$, and set the maximal iteration number as 5. In Case II, the optimality gap for MILP solver is set as 0.1%.

TABLE VIII
SOLVING RESULTS OF THE LPP-MILP MODEL BY THE MULTI-STEP SOLUTION PROCEDURE UNDER DIFFERENT ITERATION NUMBER ($\Lambda =10$) FOR THE NETWORK RECONFIGURATION SCHEME (CASE II)

| No. | Iteration number | Accumulated modeling and solving time (s) | Objective function value (MW) | $E_p^m$ (%) | $E_q^m$ (%) |
|---|---|---|---|---|---|
| 1 | 0 | 1971.1586 | 1.3201 | 452.580339 | 1902.282612 |
| 2 | 1 | 1988.3403 | 1.2820 | 0.927223 | 187.925909 |
| 3 | 2 | 2004.0396 | 1.2811 | 0.043997 | 0.673316 |
| 4 | 3 | 2020.1063 | 1.2810 | 0.008814 | 0.014142 |

TABLE IX
SOLVING RESULTS OF THE LPP-MILP MODEL BY THE MULTI-STEP SOLUTION PROCEDURE UNDER DIFFERENT ITERATION NUMBER ($\Lambda =10$) FOR THE SERVICE RESTORATION OPTIMIZATION SCHEME (CASE II)

| No. | Iteration number | Accumulated modeling and solving time (s) | Objective function value (MW) | $E_p^m$ (%) | $E_q^m$ (%) |
|---|---|---|---|---|---|
| 1 | 0 | 1926.5192 | 12.5829 | 2117.210631 | 1181.621248 |
| 2 | 1 | 1977.5440 | 12.5930 | 80.944317 | 195.343292 |
| 3 | 2 | 2019.5854 | 12.5930 | 0.209120 | 0.490495 |
| 4 | 3 | 2061.0605 | 12.5930 | 0.015434 | 0.025707 |

For the network reconfiguration scheme and the service restoration optimization scheme of the 1066-bus EDS, the detailed solving results of the LPP-MILP models by the multi-step solution procedure are shown in Table VIII and IX respectively. For both network reconfiguration and the service restoration optimization schemes, after 3 iteration processes, the error indices of the solving results are effectively reduced within the error thresholds, and the multi-step solution procedures are completed.

TABLE X
RESTORED FEEDERS, BUSES AND SELF-HEALING ISLANDS OF THE 1066-BUS EDS BY THE MULTI-STEP SOLUTION PROCEDURE UNDER DIFFERENT ITERATION NUMBER ($\Lambda =10$) FOR THE SERVICE RESTORATION OPTIMIZATION SCHEME (CASE II)

| No. | Iteration number | The number restored feeders | The number of restored buses | The number of self-healing islands |
|---|---|---|---|---|
| 1 | 0 | 317 | 324 | 6 |
| 2 | 1 | 286 | 291 | 4 |
| 3 | 2 | 278 | 283 | 4 |
| 4 | 3 | 279 | 284 | 4 |

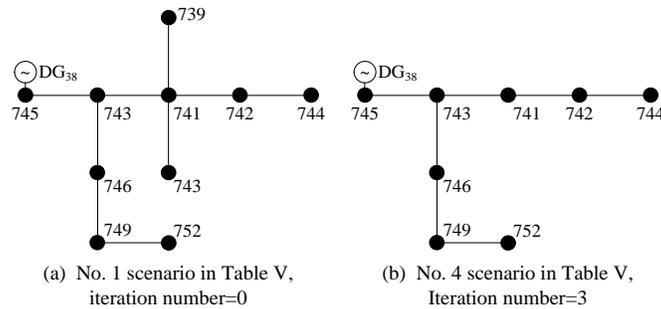

(a) No. 1 scenario in Table V, iteration number=0   (b) No. 4 scenario in Table V, Iteration number=3

Fig. 8. The topology of the self-healing island served by $DG_{38}$ (partial network of the service restoration optimization scheme of the 1066-bus EDS) under different scenarios (No.1 and No. 4) in Table IX

Take the solving results of the service restoration optimization scheme as example. The number of restored feeders, buses and self-healing islands of the 1066-bus EDS by the multi-step solution procedure are listed in Table X. As the iteration number increases, the number of restored feeders and buses gradually converges to the optimal results with progressively higher accuracy. Limited by the article space, Fig. 8 visualizes the topology of a typical self-healing island served by $DG_{38}$ (partial network of the service restoration optimization scheme of the 1066-bus EDS) under different scenarios (No.1 and No. 2) in Table IX. From Table X and Fig. 8, the solution results with relatively high approximation errors may lead to impractical or sub-optimal solving results of the service restoration optimization scheme.

The case studies are tested on a 3.2-GHz PC. With more powerful computing system, the modeling and solving time of the LPP-MILP model by the multi-step solution procedure for large EDS like the 1066-bus test system can be reduced and meet the application needs, which further demonstrate the effectiveness of the proposed multi-step solution procedure for the LPP problem in EDS.

## 6 CONCLUSIONS

In this paper, a multi-step PWL approximation based solution for LPP problem in EDS is proposed. In the multi-step solution procedure for the LPP-MILP model, the dynamic renewal strategy for the key parameters of the PWL approximation functions can efficiently reduce the approximation errors in the linearization of the network PF constraints, while the renewal logic designed in the renewal strategy can maintain the robustness of the multi-step solution procedure. And with the multi-step solution procedure, the modeling and solving time of the LPP-MILP model can be significantly decreased, which ensure the applicability of the LPP optimization scheme.

In the case studies, the proposed multi-step solution procedure for the LPP problem is tested using two real electrical distribution systems. And for each EDS test case, the two main application schemes for the LPP problem (i.e. network reconfiguration and service restoration) are tested respectively. The effectiveness of the proposed method are demonstrated via comparisons of solving results in the case studies. To stress the issue of application in relative large EDS, the proposed multi-step solution procedure for the LPP problem is further tested on a 1066-bus real distribution system in Shandong province of China.

Future work in this area will include the improvement of application of the linearized network PF constraints in multi-time scale optimization schemes. On the other hand, an unbalanced distribution system model will be considered.